\documentclass[12pt,psamsfonts,leqno,oneside,letterpaper]{amsart}
\usepackage[text={6.5truein,9truein},left=1truein,top=1truein]{geometry}
\usepackage{amssymb,amsmath,amscd,graphicx,enumerate}
\usepackage[colorlinks,linkcolor=blue,citecolor=blue,pdfstartview=FitH]{hyperref}
\usepackage{amsrefs}
\usepackage{url}
\input xy
\xyoption{all}
\SelectTips{cm}{12}

\usepackage{color}

\theoremstyle{definition}
\newtheorem{para}{}[section]

\newtheorem{remark}[para]{Remark}
\newtheorem{remarks}[para]{Remarks}
\newtheorem{notation}[para]{Notation}
\newtheorem{convention}[para]{Convention}
\newtheorem{definition}[para]{Definition}
\newtheorem{definitions}[para]{Definitions}

\theoremstyle{plain}
\newtheorem{theorem}[para]{Theorem}
\newtheorem{lemma}[para]{Lemma}
\newtheorem{proposition}[para]{Proposition}
\newtheorem{corollary}[para]{Corollary}

\newtheorem{claim}[equation]{}

\numberwithin{equation}{para}
\numberwithin{figure}{section}

\newcommand\Number{\begin{para}}
\newcommand\EndNumber{\end{para}}
\newcommand\Claim{\begin{claim}}
\newcommand\EndClaim{\end{claim}}
\newcommand\Definition{\begin{definition}}
\newcommand\EndDefinition{\end{definition}}
\newcommand\Definitions{\begin{definitions}}
\newcommand\EndDefinitions{\end{definitions}}
\newcommand\Theorem{\begin{theorem}}
\newcommand\EndTheorem{\end{theorem}}
\newcommand\Remark{\begin{remark}}
\newcommand\EndRemark{\end{remark}}
\newcommand\Remarks{\begin{remarks}}
\newcommand\EndRemarks{\end{remarks}}
\newcommand\Convention{\begin{convention}}
\newcommand\EndConvention{\end{convention}}
\newcommand\Notation{\begin{notation}}
\newcommand\EndNotation{\end{notation}}
\newcommand\Lemma{\begin{lemma}}
\newcommand\EndLemma{\end{lemma}}
\newcommand\Proposition{\begin{proposition}}
\newcommand\EndProposition{\end{proposition}}
\newcommand\Corollary{\begin{corollary}}
\newcommand\EndCorollary{\end{corollary}}
\newcommand\Proof{\begin{proof}}
\newcommand\EndProof{\end{proof}}
\newcommand\Equation{\begin{equation}}
\newcommand\EndEquation{\end{equation}}

\newcommand\Conclusions{\begin{enumerate}[(1)]}
\newcommand\EndConclusions{\end{enumerate}}

\newcommand\Properties{\begin{enumerate}[(1)]}
\newcommand\EndProperties{\end{enumerate}}
\newcommand\Conditions{\begin{enumerate}[(i)]}
\newcommand\EndConditions{\end{enumerate}}
\newcommand\Alternatives{\begin{enumerate}[(i)]}
\newcommand\EndAlternatives{\end{enumerate}}

\newcommand\Bullets{\begin{itemize}}
\newcommand\EndBullets{\end{itemize}}

\newdimen\caseindent
\global\caseindent=60pt

\newcommand\Case[1]{\par\noindent\hangindent\caseindent
     \hbox to \caseindent{\hskip .5\caseindent minus .5\caseindent
     #1\enspace\hfill}\ignorespaces}

\newcommand\boibsectionone{2}
\newcommand\goodtowerdef{8.4}
\newcommand\topeleven{8.13}

\newcommand\geomeleven{9.6}
\newcommand\threefreevolume{9.3}
\newcommand\firstAST{9.4}
\newcommand\moosday{2.11}
\newcommand\newimprovedoldpropthree{4.4}
\newcommand\towerproposition{8.11}
\newcommand\simpletower{8.12}
\newcommand\tworminusfour{8.5}
\newcommand\card{\mathop{\#}}
\newcommand\Size{size}

\newcommand\rk{\mathop{{\rm rk}_2}}

\newcommand\vol{\mathop{\rm Vol}}

\newcommand\genus{\mathop{\rm genus}}

\newcommand\Z{{\mathbb Z}}
\newcommand\ZZ{{\mathbb Z}}

\newcommand\cala{{\mathcal A}}

\newcommand\calb{{\mathcal B}}
\newcommand\calc{{\mathcal C}}

\newcommand\caln{{\mathcal N}}
\newcommand\calp{{\mathcal P}}

\newcommand\cals{{\mathcal S}}
\newcommand\calt{{\mathcal T}}

\newcommand\calu{{\mathcal U}}

\newcommand\calw{{\mathcal W}}
\newcommand\calx{{\mathcal X}}
\newcommand\caly{{\mathcal Y}}
\newcommand\calz{{\mathcal Z}}
\newcommand\inter{\mathop{\rm int}}

\newcommand\barchi{\bar\chi}
\newcommand\chibar{\bar\chi}

\newcommand{\tN}{\widetilde N}

\begin{document}

\author{Marc Culler}
\address{Department of Mathematics, Statistics, and Computer Science (M/C 249)\\
University of Illinois at Chicago\\
851 S. Morgan St.\\
Chicago, IL 60607-7045}
\email{culler@math.uic.edu}
\thanks{Partially supported by NSF {grants  DMS-0204142 and DMS-0504975}}

\author{Peter B. Shalen}
\address{Department of Mathematics, Statistics, and Computer Science (M/C 249)\\
University of Illinois at Chicago\\
851 S. Morgan St.\\
Chicago, IL 60607-7045}
\email{shalen@math.uic.edu}

\title{Singular surfaces, mod 2 homology,  and hyperbolic volume, II}

\begin{abstract}
If $g$ is an integer $\ge2$, and $M$ is a closed simple $3$-manifold
such that $\pi_1(M)$ has a subgroup isomorphic to a genus-$g$ surface
group and $\dim_{\Z_2}H_1( M;\Z_2) \ge \max(3g-1,6)$, we show that $M$
contains a closed, incompressible surface of genus at most $g$. As an
application we show that if $M$ is a closed orientable hyperbolic
$3$-manifold such that $\vol M \le3.08$, then $\dim_{\Z_2}H_1( M;\Z_2)
\le 5$.
\end{abstract}

\maketitle

\section{Introduction}

This paper is a sequel to \cite{last}. As in \cite{last}, we
write $\rk V$ for the dimension of a $\Z_2$-vector space $V$, and set
$\rk X=\rk H_1(X;\Z_2)$ when $X$ is a space of the homotopy type of a
finite CW-complex. As in \cite{last}, we say that an orientable
$3$-manifold $M$ is {\it simple} if $M$ is compact, connected,
orientable, irreducible and boundary-irreducible, no subgroup of
$\pi_1(M)$ is isomorphic to $\Z\times\Z$, and $M$ is not a
closed manifold with finite fundamental group.

We shall establish the following topological result,
which is a refinement of Theorem \topeleven\ of \cite{last}.

\Theorem\label{top 6}
Let $g$ be an integer $\ge2$.
Let $M$ be a closed simple $3$-manifold such that $\rk M \ge \max(3g-1,6)$ and
$\pi_1(M)$ has a
subgroup isomorphic to a genus-$g$ surface group. Then $M$ contains a
closed, incompressible surface of genus at most
$g$.
\EndTheorem

Like \cite{last}*{Theorem \topeleven}, this result may be regarded
as a partial analogue of Dehn's lemma for $\pi_1$-injective genus-$g$
surfaces. The difference between the two theorems is that
the hypothesis $\rk M \ge \max(3g-1,6)$ assumed in Theorem \ref{top 6} 
is strictly weaker than the corresponding hypothesis in
\cite{last}*{Theorem \topeleven}, namely that $\rk M \ge 4g-1$.

For the case $g=2$, Theorem \ref{top 6} is almost sharp: in Section
\ref{example section} we construct examples of simple $3$-manifolds
$M$ with $\rk M=4$ such that 
$\pi_1(M)$ has a
subgroup isomorphic to a genus-$2$ surface group, but $M$ contains no
closed, incompressible surface whatever.

As an application of Theorem \ref{top 6} we shall prove the following
theorem relating volume to homology for closed hyperbolic
$3$-manifolds.

\Theorem\label{geom 6}
Let $M$ be a closed orientable hyperbolic $3$-manifold such that
$\vol M \le3.08$. Then $\rk M \le 5$.
\EndTheorem

Theorem \ref{geom 6} is a refinement of Theorem \geomeleven\ of
\cite{last}, and will be deduced from Theorem \ref{top 6} in the same
way that \cite{last}*{Theorem \geomeleven} was deduced from
\cite{last}*{Theorem \topeleven}.  

In \cite{4-free}, by combining Theorem \ref{top 6} with new geometric
results, we will prove that if $M$ is a closed orientable hyperbolic
$3$-manifold such that $\vol M \le3.44$, then $\rk M \le 7$.
Further applications of Theorem \ref{top 6} to the study of volume and
homology will be given in \cite{cds}.

The arguments in this paper draw heavily on results from
\cite{last}. The improvements that we obtain here depend on a much
deeper study of books of $I$-bundles (see \cite{last}*{Section
  \boibsectionone}) in closed $3$-manifolds than the one made in
\cite{last}. For all $g\ge2$ this involves new topological
ingredients. For $g>2$ it also involves a surprising application of
Fisher's inequality from combinatorics.

Before describing the new ingredients in the proof of Theorem \ref{top
  6} we shall briefly review the proof of \cite{last}*{Theorem
  \topeleven} and explain the role that books of $I$-bundles play in
it.  The proof uses a tower of two-sheeted covers analogous to the one
used by Shapiro and Whitehead in their proof of Dehn's lemma
\cite{ShapiroWhitehead}. The homological hypothesis allows one to
construct a good tower (in the sense of \cite{last}*{Definition
  \goodtowerdef})
$${\mathcal T}=(M_0,N_0,p_1,M_1,N_1,p_2,\ldots,p_n,M_n,N_n),$$ with
base $M_0$ homeomorphic to $M$ and with some height $n\ge0$, such that
{$N_n$} contains a connected (non-empty) closed incompressible surface
$F$ of genus $\le g$. (Here $N_j$ is a submanifold of $M_j$ for
$j=0,\ldots,n$ and $p_j:M_j\to N_{j-1}$ is a two-sheeted covering map
for $j=1,\ldots,n$.) The key step is to show, for a given $j>0$, that
if $N_j$ contains a connected closed incompressible surface $F$ of
genus $\le g$, then $N_{j-1}$ contains such a surface as well. Certain
books of $I$-bundles arise as obstructions to carrying out this step.
Specifically, the arguments of \cite{last} show that this step can be
carried out unless $N_{j-1}$ is a closed manifold that contains a
submanifold of the form $W=|\calw|$, where $\calw$ is a book of
$I$-bundles, $\chi(W)\ge 2-2g$, and the inclusion homomorphism
$H_1(W;\Z_2)\to H_1(N_{j-1};\Z_2)$ is surjective.  This situation is
then ruled out by estimating ranks of homology groups.  Under the
hypothesis of \cite{last}*{Theorem \topeleven}, one can show that when
$N_{j-1}$ is closed we have $\rk N_{j-1}\ge4g-2$, while
\cite{last}*{Lemma \moosday} implies that $\rk W\le4g-3$ when
$\chi(W)\ge2-2g$.  Thus the induced map on homology cannot be
surjective.

Under the weaker hypothesis of Theorem \ref{top 6} one obtains only a
lower bound of $\max(3g-2,5)$ for $\rk N_{j-1}$ when $ N_{j-1}$ is
closed.  So the homological condition given by Lemma \moosday\ of
\cite{last} does not suffice to overcome the obstruction. Instead,
the strategy for carrying out the key step is to first attempt
to construct the required incompressible surface by compressing
the boundary of a carefully chosen submanifold of $W$.

It is easy to choose the book of $I$-bundles $\calw$ defining $W$ so
that each of its pages has Euler characteristic $-1$. In this case one
can find a sub-book $\calw_0$ of $\calw$ such that $W_0=|\calw_0|$ has
exactly half the Euler characteristic of $W$. Using classical
$3$-manifold techniques one can then show that either (a) the
inclusion homomorphism $\iota_\sharp :\pi_1(W_0)\to\pi_1(N_{j-1})$ has image of
rank at most $g$, or (b) $\iota_\sharp$  is surjective, or (c) a connected
incompressible surface can be obtained from $\partial W$ by doing
ambient surgeries in $N_{j-1}$ and selecting a component. One can use
Lemma \moosday\ of \cite{last} to show that alternative (b)
contradicts the lower bound for $\rk N_{j-1}$. If alternative (c)
holds, one has an incompressible surface of genus less than $g$, which
is all that the tower argument requires. If (a) holds, a relative
version of the proof of \cite{last}*{Lemma \moosday} gives an upper
bound of $3g-2$ for $\rk N_{j-1}$;  this contradicts our condition
$\rk N_{j-1} \ge \max(3g-2,5)$ unless $g>2$ and $\rk N_{j-1}=3g-2$.

To deal with the latter situation we must exercise even more care in
choosing the sub-book $\calw_0$. It turns out  (see Lemma \ref{we'll
  all go to the seashore})  that when  $g>2$  one can choose
$\calw_0$ in such a way that  $H_2(W_0;\Z_2)\ne0$.  In particular it then
follows that $W_0$ is not a handlebody, and the classical $3$-manifold
argument mentioned above can be modified to show that either (b) or
(c) holds, or else (a$'$) the image of $\iota_\sharp$ has rank at most
$g-1$. One can then improve the upper bound for $\rk N_{j-1}$ to
$3g-2$ and obtain the required contradiction.

Making the right choice for $\calw_0$ in this case requires both a
detailed study of the homology of books of $I$-bundles and an
interesting result, Proposition \ref{quality control}, about
finite-dimensional vector spaces over $\Z_2$. It is in the proof
of Proposition \ref{quality control} that we need to apply
Fisher's inequality.

Section \ref{compressing section} contains the classical $3$-manifold
arguments that we mentioned in the outline above, and Section
\ref{algebra section}  is devoted to the proof of Proposition
\ref{quality control}. In Section \ref{kalamazoo to timbuctu}
these ingredients are combined with some observations about homology
of books of $I$-bundles to carry out the main step, sketched above, in
the proof of Theorem \ref{top 6}; the proof of the theorem itself
appears in Section \ref{illegitimis non carborundum}. Section
\ref{example section} is devoted to constructing the examples,
referred to above, that show
that the theorem is almost sharp.

In Section \ref{non-fibroid section} we establish a stronger version
of Theorem \ref{top 6}, Proposition \ref{stealing that extra bow}, 
which is particularly well-adapted to the applications to volume
estimates, including the proof of Theorem \ref{geom 6} and the
application  in the forthcoming paper \cite{4-free}. The proof of
Theorem \ref{geom 6}  is given in Section \ref{volume section}.

In general we will  use all of the conventions that were used in
\cite{last}. In particular, in addition to the notations $\rk V$ and
$\rk X$, and the definition of a simple manifold, which were explained
above, we shall set $\chibar(X)=-\chi(X)$ when $X$ is a space of the
homotopy type of a finite CW-complex (and $\chi(X)$ as usual denotes
its Euler characteristic). {\it Connected} spaces are understood to be
in particular non-empty, and {\it irreducible} $3$-manifolds are
understood to be in particular connected.

The cardinality of any finite set $S$ will be denoted 
by $\card S$.

We are grateful to Ian Agol for many valuable discussions, and of
course for his crucial contribution to \cite{last}. We are also
grateful to Dhruv Mubayi for telling us about Fisher's inequality.
Finally, we thank the anonymous referee for an impressively prompt and
thorough job of reviewing the paper and for asking a question
which led to our discovery of the material in Section \ref{example
  section}.

\section{Compressing submanifolds}\label{compressing section}

Recall that  a {\it
compressing disk} for a closed surface $F$ in the interior of a
$3$-manifold  $ M$ is defined to be a disk $D\subset M$ such that
$D\cap F=\partial D$, and such that $\partial D$ does not
bound a disk in $F$.

\Lemma\label{wicopee}Let $M$ be a compact, connected, orientable,
simple $3$-manifold, and let $T$ be a compressible torus in $\inter M$.
Then either $T$ bounds a solid torus in $\inter M$, or $T$ is contained in a
ball in $\inter M$.
\EndLemma

\Proof Since $M$ is simple, $T$ is compressible. Fix a compressing
disk $D$ for $X$. Let $E\subset\inter M$ be a ball containing $D$,
such that $A\doteq E\cap T\subset \partial E$, and such that $A$
is a regular neighborhood of $\partial D$ in $T$. Then
$\overline{(\partial E)-A}$ has two components $D_1$ and $D_2$, both
of which are disks, and $D_1\cup A\cup D_2$ is a sphere, which must
bound a ball $B\subset M$. By connectedness we have either $E\cap
B=D_1\cup D_2$ or $E\subset B$. The first alternative implies that
$E\cup B\supset T$ is a solid torus, and the second implies that
$T\subset B$.
\EndProof

\Definitions\label{all chalk, no traction}
Let $M$ be a compact, connected, orientable, irreducible $3$-manifold.
We shall denote by $\calx_M$ the set of all compact, connected
$3$-submanifolds $X$ of $\inter M$ such that \Conditions
\item no component of $\partial X$ is a $2$-sphere, and 
\item $X$ does not carry $\pi_1(M)$, i.e. the inclusion
homomorphism $\pi_1(X)\to\pi_1(M)$ is not surjective.
\EndConditions

For any $X\in\calx_M$, since $\partial X$ has no $2$-sphere
components, we have $\chibar(X)\ge0$. We let $t(X)$ denote the number
of components of $\partial X$ that are tori, and we set
$$k(X)=t(X)+3\chibar(X)\ge0.$$

For any $X\in\calx_M$ we denote by $r(X)$ the rank of the image of the
inclusion homomorphism $\pi_1(X)\to\pi_1(M)$. We set
$$i(X)=\chibar(X)-r(X)\in\ZZ.$$

If $X\in\calx_M$ is given, we define a {\it compressing disk} for
$X$ to be a compressing disk for $\partial X$. We shall say that $D$
is {\it internal} or {\it external} according to whether $D\subset X$ or
$D\cap X=\partial D$.   We shall say that an internal compressing disk 
 $D$ is {\it separating} if $X-D$ is connected, and
{\it non-separating} otherwise. 
\EndDefinitions

\Lemma\label{gentle jane was as good as gold} 
Let $M$ be a compact, connected, orientable, irreducible $3$-manifold,
and let $X\in\calx_M$ be given, Suppose that every component of
$\partial X$ has genus strictly greater than $1$, and that $X$ has an
internal compressing disk.  Then there is an element $X'$ of $\calx_M$
such that \Conclusions
\item\label{taffy was a welshman} $X'\subset X$,
\item\label{mama puts it in my milk} every component of $M-X'$ contains
  at least one component of $M-X$,
\item\label{or fostered a passion for alcohol} $\chibar(X')\le\chibar(X)-1$,
\item\label{taffy was a thief} $k(X')<k(X)$,
\item\label{taffy came to my house} $i(X')\le\max( i(X),( i(X)-1)/2)$, and
\item\label{to steal a bit of beef} if $X$ is not a handlebody then
  $X'$ is not a handlebody.
\EndConclusions
\EndLemma

\Proof
If $X$ has a non-separating internal compressing disk we fix such a
disk and denote it by $D$.  If every internal compressing disk for $X$
is separating we let $D$ denote an arbitrarily chosen  internal
compressing disk for $X$. In either case we set $\gamma=\partial D$,
and denote by $F$ the component of $\partial X$ that contains
$\gamma$.

We let $E$ denote a regular neighborhood of $D$ in $X$. The manifold
$Z=\overline{X-E}$ has at most two components. Each component of
$\partial Z$ is either a component of $\partial X$, or a component of
the surface obtained from $F$ by surgery on the simple closed curve
$\gamma$, which is homotopically non-trivial in $F$. Since every
component of $\partial X$ has genus strictly greater than $1$, it
follows that no component of $\partial Z$ is a $2$-sphere, and that
$\partial Z$ has at most two torus components. Since
$Z\subset X$ and $X\in\calx_M$, no component of $Z$ carries
$\pi_1(M)$. Hence each component of $Z$ belongs to $\calx_M$.

We observe that
\Equation\label{i had a little nut tree}
Z\subset X,
\EndEquation
that
\Equation\label{i fooled mama and put it in her tea}
\text{every component of $M-Z$ contains
  at least one component of $M-X$,}
\EndEquation
and that
\Equation\label{but nothing would it bear}
\chibar(Z)=\chibar(X)-1.
\EndEquation

Since $\partial Z$ has at most two torus components, we have 
\Equation\label{but a silver nutmeg}t(Y)\le2\EndEquation
for every component $Y$ of $Z$.

We claim:
\Claim\label{give that man a kewpie doll} There is a component $X'$ of
$Z$ such that
$i(X')\le\max( i(X),( i(X)-1)/2)$.
\EndClaim

We first prove \ref{give that man a kewpie doll} in the case where $Z$
is connected. We shall show that $i(Z)\le i(X)$, which implies
\ref{give that man a kewpie doll} in this case.
We fix a base point
$\star\in Z$ and let $G,G'\le\pi_1(M,\star)$ denote the respective
images of $\pi_1(X,\star)$ and $\pi_1(Z,\star)$ under inclusion. Then
$G$ is generated by $G'$ and $\alpha$, where $\alpha\in\pi_1(M,\star)$
is the homotopy class of a loop in $X$ that crosses $D$ in a single
point. Hence $r(X)\le r(Z)+1$. In view of (\ref{but nothing would it
  bear}), it follows that
$i(X)\ge i(Z)$.

We next prove \ref{give that man a kewpie doll} in the case where $Z$
is disconnected. Let $Y_1$ and $Y_2$ denote the components of $Z$.
We fix a base point $\star\in D$ and let $G,G_i'\le\pi_1(M,\star)$
denote the respective images of $\pi_1(X,\star)$ and
$\pi_1(Y_i,\star)$ under inclusion. Then $G$ is generated by $G_1'$
and $G_2'$, so that $r(X)\le r(Y_1)+r(Y_2)$. 
It follows that 
$$\begin{aligned}i(Y_1)+i(Y_2)&=(\chibar(Y_1)+\chibar(Y_2))-(r(Y_1)+r(Y_2))\\
&=(\chibar(Y_1)+\chibar(Y_2))-(r(Y_1)+r(Y_2))\\
&\le\chibar(Z)-r(X)\end{aligned}$$
which in view of (\ref{but nothing would it
  bear}) gives
$$i(Y_1)+i(Y_2)\le i(X)-1.$$
Hence for some $j\in\{1,2\}$ we have
$$i(Y_j)\le \frac{i(X)-1}2.$$
If we set $X'=Y_j$ with this choice of $j$, then \ref{give that man a
  kewpie doll} follows in this case.

Now let $X'$ denote the component of $Z$ given by \ref{give that man a
  kewpie doll}. Thus Conclusion (\ref{taffy came to my house}) of the
lemma holds with this choice of $X'$. In view of \ref{i had a little
  nut tree}, Conclusion (\ref{taffy was a welshman}) holds as well.

It follows from \ref{but nothing would it bear} that
$\chibar(X)-1=\sum_Y\chibar(Y)$, where $Y$ ranges over the components
of $Z$. Since each component of $Z$ belongs to $\calx_M$, we have
$\chibar(Y)\ge0$ for each component $Y$ of $Z$.  Hence 
$\chibar(Y)\le\chibar(X)-1$ for each component $Y$ of $Z$. This,
together with (\ref{but a silver nutmeg}), implies that
$k(Y)<k(X)$ for each component $Y$ of $Z$.
In particular, Conclusions (\ref{or fostered a passion for alcohol}) and (\ref{taffy was a thief}) hold with
our choice of $X'$.

Since $X'$ is a component of $Z$, every component of $M-X'$ contains
at least one component of $M-Z$. Combining this observation with
(\ref{i fooled mama and put it in her tea}) we obtain Conclusion
\ref{mama puts it in my milk}.

It remains to prove Conclusion (\ref{to steal a bit of beef}). We
shall assume that $X'$ is a handlebody and deduce that $X$ is a
handlebody. If $Z$ is connected, so that $X'=Z$, then $X$ is the union
of the handlebody $Z$ and the ball $E$, and $Z\cap E$ is the union of
two disjoint disks. Hence $X$ is a handlebody. Now suppose that $Z$ is
disconnected, i.e. that $X-D$ is disconnected. Since the handlebody
$X'$ is an element of $\calx_M$, it must have strictly positive genus.
Hence there is a disk $D'\subset X'$ such that $X'-D'$ is connected.
After modifying $D'$ by an ambient isotopy in $X'$ we may assume that
$D'$ is disjoint from the disk $X'\cap E\subset\partial X'$. Then $D'$
is an internal compressing disk for $X$, and $X-D'$ is connected. But
in this case the choice of $D$ guarantees that $X-D$ is
connected, and we have a contradiction.
\EndProof

\Lemma\label{she always did as she was told}
Let $M$ be a compact, connected, orientable,
irreducible $3$-manifold, and let $X\in\calx_M$ be given, Suppose that
every component 
of $\partial X$ has genus strictly greater
than $1$, and that $X$ has an external compressing disk.  Then there is an element $X'$ of $\calx_M$ such that
\Conclusions
\item\label{or vivisected her last new doll}  $\chibar(X')=\chibar(X)-1$,
\item\label{she never spoke when her mouth was full}
$k(X')<k(X)$, and
\item\label{or caught bluebottles their legs to pull} $i(X')=
  i(X)-1$.
\EndConclusions
\EndLemma

\Proof We fix an external compressing disk $D$ for $X$, we set
$\gamma=\partial D$, and we let $E$ denote a regular neighborhood of
$D$ in $\overline{M-X}$. We set $X'=X\cup E$. Note that the inclusion
homomorphism $\pi_1(X)\to\pi_1(X')$ is surjective. Hence:
\Claim\label{or spilt plum jam on her nice new frock} For any base
point $\star\in X$, the inclusion homomorphisms
$\pi_1(X,\star)\to\pi_1(M,\star)$ and
$\pi_1(X',\star)\to\pi_1(M,\star)$ have the same image.  \EndClaim
Since $X\in\calx_M$, it follows from \ref{or spilt plum jam on her
nice new frock} that $X'$ does not carry $\pi_1(M)$. On the other
hand, each component of $\partial X'$ is either a component of
$\partial X$, or a component of the surface obtained from $F$ by
surgery on the simple closed curve $\gamma$, which is homotopically
non-trivial in $F$. Since every component of $\partial X$ has genus
strictly greater than $1$, it follows that no component of $\partial
X'$ is a $2$-sphere, and that at most two of its components are tori.
Hence $X'\in\calx_M$, and \Equation\label{and a golden
pear}t(X')\le2.\EndEquation

With this definition of $X'$, it is clear that  Conclusion
(\ref{or vivisected her last new doll}) of the lemma holds. With
(\ref{and a golden pear}), this implies Conclusion (\ref{she never
spoke when her mouth was full}).    On the other hand, by
\ref{or spilt plum jam on her nice new frock} we have
\Equation\label{or put white mice in the eight-day clock} r(X')=r(X).
\EndEquation Combining (\ref{or put white mice in the eight-day
clock}) with Conclusion (\ref{she never spoke when her mouth was
full}), we immediately obtain Conclusion (\ref{or caught bluebottles
their legs to pull}).  \EndProof

\Definition\label{who the hell said anything about christmas}
Let $g\ge2$ be an integer and let $M$ be a closed, orientable,
irreducible $3$-manifold. We shall say that $M$ is {\it $g$-small} if
$M$ contains no separating, closed, incompressible surface of genus
$g$, and contains no closed incompressible surface of genus $<g$.
\EndDefinition

\Lemma\label{and who the hell did}
Let $c$ be a positive integer, let $M$ be a compact, connected,
orientable, irreducible $3$-manifold which is $ (c+1)$-small, and let
$X$ be an element of $\calx_M$ such that $\chibar(X)\le c$.  Then every
component of $\partial X$ is compressible in $M$.
\EndLemma

\Proof The hypothesis that $\chibar(X)\le c$ implies that every
component of $\partial X$ has genus at most $c+1$, and that if
$\partial X$ is disconnected then each of its components has genus at
most $c$. In particular, every component of $\partial X$ is either a
separating surface of genus $c+1$, or a surface of genus at most $c$.
Since $M$ is $(c+1)$-small it follows that every component of
$\partial X$ is compressible in $M$.
\EndProof

\Proposition\label{new but not newest old prop 2}
Let $M$ be a compact, connected, orientable, irreducible $3$-manifold,
and let $Y$ be an element of $\calx_M$.  Set $c=\chibar(Y)$,
and assume that $M$ is $ (c+1)$-small.  Then $i(Y)\ge-1$.
\EndProposition

\Proof
Suppose that $i(Y)\le-2$. Let $\calx^*_M\subset\calx_M$ denote the set
of all $X\in\calx_M$ such that
\Conditions
\item\label{and when she grew up}$\chibar(X)\le c$ and
\item\label{she was given in marriage}$i(X)\le-2$.
\EndConditions
Then $Y\in\calx^*_M$ and so $\calx^*_M\ne\emptyset$. Let us choose an
element $X\in\calx^*_M$ such that $k(X)\le k(X')$ for every
$X'\in\calx^*_M$.

Since $X$ belongs to $\calx_M$, it cannot carry $\pi_1(M)$; in
particular, $X\ne M$, and so $\partial X\ne\emptyset$. Since
$X\in\calx^*_M$, we have $\chibar(X)\le c$. It therefore follows from
Lemma \ref{and who the hell did} that every component of $\partial X$
is compressible in $M$. In particular $X$ has either an internal or an
external compressing disk.

We first consider the case in which $X$ has an internal
compressing disk, and every component of $\partial X$ has genus
$>1$. In this  case, Lemma \ref{gentle jane was as good as gold} gives
an element $X'$ of $\calx_M$ such that
$\chibar(X')\le\chibar(X)-1$,
$k(X')<k(X)$, and $i(X')\le\max( i(X),( i(X)-1)/2)$. Since
$\chibar(X)\le c$ and $i(X)\le-2$, it follows that 
$\chibar(X')\le c-1$ and that $i(X')<-1$. In particular,
$X'\in\calx^*_M$. Since $k(X')<k(X)$, this contradicts our choice of
$\calx$. 

We next turn to the case in which $X$ has an external
compressing disk, and every component of $\partial X$ has genus
$>1$. In this  case, Lemma 
\ref{she always did as she was told}
gives
an element $X'$ of $\calx_M$ such that
$\chibar(X')\le\chibar(X)-1$,
$k(X')<k(X)$, and $i(X')=i(X)-1$. Since
$\chibar(X)\le c$ and $i(X)\le-2$, it again follows that 
$\chibar(X')\le c-1$ and that $i(X')<-1$. Again it follows that
$X'\in\calx^*_M$, and since $k(X')<k(X)$,  our choice of
$\calx$ is contradicted.

There remains the case in which some component $T$ of $\partial X$ is a
torus. According to Lemma \ref{wicopee},  $T$ is the boundary of a
compact submanifold $W$ of $\inter M$ such that either (a) $W$ is a solid
torus, or (b) $W$ is contained in a
ball in $\inter M$. We must have either $X\subset W$ or $X\cap W=T$.

Either of the alternatives (a) or (b) implies that the image of
$\pi_1(W)$ under the inclusion to $\pi_1(M)$ is at most cyclic. Hence
if $X\subset W$ then $r(X)\le1$, and hence $i(X)\ge-1$. This is a
contradiction since $X\in\calx^*_M$.

If $X\cap W=T$, we set $X'=X\cup W$. Then $\partial X'=(\partial
X)-T$. In particular $\partial X'$ has no sphere components. If
$\star$ is a base point in $X$, either of the alternatives (a) or (b)
implies that $\pi_1(X,\star)$ and $\pi_1(X',\star)$ have the same
image under the inclusion to $\pi_1(M,\star)$. It follows that $X'$
does not carry $\pi_1(M)$, so that $X\in\calx_M$. It also follows
that $r(X')=r(X)$. But since $\partial X'=(\partial X)-T$, we have
$\chibar(X')=\chibar(X)$ and $t(X')=t(X)-1$. We now deduce that 
$\chibar(X')\le c$ and $i(X')=i(X)\le-2$, so that $X'\in\calx^*_M$; and
that $k(X')=k(X)-1$. This contradicts our choice of $X$.
\EndProof

\Proposition\label{newest old prop 2}Let $M$ be a compact, connected,
orientable, irreducible $3$-manifold, and let $Y$ be an element of
$\calx_M$. Assume that $Y$ is not a handlebody and that no component
of $\partial Y$ is a torus.  Set $c=\chibar(Y)$, and assume that $M$ is $
(c+1)$-small.  Then $i(Y)\ge0$.  \EndProposition

\Proof
We define a {\it $Y$-special submanifold} of $M$ to be a compact
$3$-dimensional submanifold $W$ of $M$ such that
\Bullets
\item $\partial W$ is a torus contained in $\inter Y$,
\item $W\not\subset Y$, and
\item either $W$ is a solid torus or $W$ is contained in a ball in
  $\inter M$.
\EndBullets

We distinguish two cases.

{\bf Case I. There is no $Y$-special submanifold of $M$.}

In order to prove that in Case I we have $i(Y)\ge0$, we reason by
contradiction. Assume that $i(Y)\le-1$.
Let $\calx^{**}_M\subset\calx_M$ denote the set of
all $X\in\calx_M$ such that
\Conditions
\item\label{teasing tom was a very bad boy}$X\subset \inter Y$,
\item\label{a great big squirt was his favorite toy}every component of $M-X$ contains
  at least one component of $M-Y$,
\item\label{he punched his little sisters' heads}$X$ is not a handlebody,
\item\label{to a first-class earl}$\chibar(X)\le c$ and
\item\label{who keeps his carriage}$i(X)\le-1$.
\EndConditions
The hypotheses
and the assumption that $i(Y)\le-1$,
imply that a manifold obtained from $Y$ by removing a half-open collar about
$\partial Y$ belongs to $\calx^{**}_M$. Hence $\calx^{**}_M\ne\emptyset$. Let us choose an
element $X\in\calx^{**}_M$ such that $k(X)\le k(X')$ for every
$X'\in\calx^{**}_M$.

Since $X$ belongs to $\calx_M$, it cannot carry $\pi_1(M)$; in
particular, $X\ne M$, and so $\partial X\ne\emptyset$. Since
$X\in\calx^{**}_M$, we have $\chibar(X)\le c$.
It therefore follows from
Lemma \ref{and who the hell did} that every component of $\partial X$
is compressible in $M$.
 In particular, $X$ has either an internal or an external
compressing disk.

We first consider the subcase in which $X$ has an internal compressing
disk, and every component of $\partial X$ has genus $>1$. In this
case, there is an element $X'$ of $\calx_M$ such that Conclusions
(\ref{taffy was a welshman})--(\ref{to steal a bit of beef}) of
Lemma \ref{gentle jane was as good as gold} hold.
Since $X\in\calx^{**}_M$, Conclusions (\ref{taffy was a welshman}),
(\ref{mama puts it in my milk}),
(\ref{or fostered a passion for alcohol}),
(\ref{taffy came to my house}) and
(\ref{to steal a bit of beef}) of Lemma  \ref{gentle jane was as good
as gold} imply, respectively, that $X'\subset Y$; that every component of $M-X'$ contains
  at least one component of $M-Y$; that $\chibar(X')\le c-1$; that 
$i(X')\le-1$; and that $X'$ is not a handlebody. Hence
$X'\in\calx^{**}_M$. But (\ref{taffy was a thief}) gives
$k(X')<k(X)$, and this contradicts our choice of
$\calx$. 

We next turn to the subcase in which $X$ has an external
compressing disk, and every component of $\partial X$ has genus
$>1$. In this  case, Lemma 
\ref{she always did as she was told}
gives
an element $X'$ of $\calx_M$ such that
$\chibar(X')\le\chibar(X)-1$ and $i(X')=i(X)-1$. Let us set
$c'=\chibar(X')$. Since
$X\in\calx^{**}_M$ we have $c'\le\chibar(X)-1\le c-1$. Since by
hypothesis $ M$ contains no incompressible closed surfaces of
genus $\le c+1$, in particular it contains no 
incompressible closed surfaces of
genus $\le c'+1$. Hence the hypotheses of Proposition \ref{new but not
newest old prop 2} hold with $X$ and $c'$ in place of $Y$ and $c$. It
follows that $i(X')\ge-1$, i.e. that $i(X)\ge0$. But since
$X\in\calx^{**}_M$ we have
$i(X)\le-1$, a contradiction.

The remaining subcase of Case
I is the one in which some component $T$ of $\partial X$ is a torus.
According to Lemma \ref{wicopee}, $T$ is the boundary of a compact
submanifold $W$ of $\inter M$ such that either  $W$ is a solid
torus, or $W$ is contained in a ball in $\inter M$. Since we are
in Case I, the submanifold $W$ of $M$ cannot be $Y$-special. Hence we
must have $W\subset Y$.

Since $\partial W=T\subset\partial X$, we must have either $X\subset
W$ or $X\cap W=T$.  If $X\cap W=T$, then $\inter W$ is a component of
$M-X$. By Condition (\ref{a great big squirt was his favorite toy}) in
the definition of $\calx^{**}_M$, the set $W$ must contain a component
of $M-Y$. This is impossible since $W\subset Y$.

Now suppose that $X\subset W$. If $X$ is a proper subset of $W$ then
$W$ contains a component of $M-X$, which by Condition (\ref{a great big
squirt was his favorite toy}) in the definition of $\calx^{**}_M$ must
contain a component of $M-Y$.  Again this is impossible since
$W\subset Y$. Hence $X=W$. If $W$ is a solid torus, we have a
contradiction to Condition (\ref{he punched his little sisters' heads})
in the definition of $\calx^{**}_M$. Finally, if $X$ is contained in a
ball in $\inter M$ we have $r(X)=0$ and hence $i(X)\ge0$.
This contradicts Condition (\ref{who keeps his carriage})
in the definition of $\calx^{**}_M$. 

{\bf Case II. There is a $Y$-special submanifold of $M$.}
In this case, we fix a $Y$-special submanifold $W$ of $M$, and we set
$Y'=Y\cup W$. We also set $F=W\cap\partial Y$. The definition of a
$Y$-special manifold guarantees that $W\not\subset Y$ and hence that
$F\ne\emptyset$. But $F$ is a union of components of $\partial Y$, and
by the hypothesis of the proposition, no component of $\partial Y$ is
a torus. Hence $\chibar(F)>0$.

We have $\partial Y'=(\partial Y)-F$. In particular $\partial Y'$ has
no sphere components.

According to the definition of a $Y$-special submanifold, either $W$
is a solid torus or $W$ is contained in a ball in $\inter M$. In
either case, if $\star$ is a base point in $Y$, then $\pi_1(Y,\star)$
and $\pi_1(Y',\star)$ have the same image under the inclusion to
$\pi_1(M,\star)$. It follows that $Y'$ does not carry $\pi_1(M)$, so
that $Y\in\calx_M$. It also follows that $r(Y')=r(Y)$. But since
$\partial Y'=(\partial Y)-F$, we have
$\chibar(Y')=\chibar(Y)-\chibar(F)<\chibar(Y)$. It follows that
\Equation\label{and dropped hot ha'pennies down their backs}
i(Y')<i(Y).
\EndEquation

On the other hand, if we set $c''=\chibar(Y')<\chibar(Y)=c$, the
hypothesis of the proposition implies that $M$
contains no incompressible closed surface of genus $\le c''+1$. Hence
by Proposition \ref{new but not newest old prop 2} we have 
$i(Y')\ge-1$. In view of (\ref{and dropped hot ha'pennies down their
  backs}) it follows  that $i(Y)\ge0$.
\EndProof

\section{An algebraic result}\label{algebra section}
\def\size_#1(#2){||#2||_#1}

Suppose that $V$ is a finite-dimensional vector space over $\ZZ_2$ and
that $\calu$ is a basis of $V$. Then any element  $\alpha$ of
$V$ may be written uniquely in the form $\sum_{u\in\calu}\lambda_uu$,
with $\lambda_u\in\ZZ_2$ for each $u\in\calu$. We denote by 
$S_\calu(\alpha)$  the set of elements $u\in\calu$ such that 
$\lambda_u=1$,  and define the {\it \Size} of $\alpha$ with respect to
the basis $\calu$, denoted $\size_\calu(\alpha)$, to be 
$\card{S_\calu(\alpha)}$.   Note that
$\size_\calu(\alpha)=0$  if and only if $\alpha=0$.

The purpose of this section is to prove the following result:

\Proposition\label{quality control}
Let $m$ be an integer $\ge2$, let $\calu$ be a basis of a
$2m$-dimensional vector space $V$ over $\ZZ_2$, and suppose that $H$
is a subspace of $V$ with dimension at least $m$. Then there is an
element $\alpha$ of $H$ such that $0<\size_\calu(\alpha)\le m$.
\EndProposition

\Proof
We divide the argument into two  cases.

{\bf Case I. There is an element $\beta$ of $H$ such that
$\size_\calu(\beta)\ge m+2$.}

Set $k=\size_\calu(\beta)$, so that $m+2\le k\le 2m$.  We let $L$
denote the linear subspace of $V$ spanned by $S_\calu(\beta)$. Thus
$L$ consists of all elements $\alpha\in V$ such that $S_\calu(\alpha)\subset
S_\calu(\beta)$.

We have $\rk L=k$, and so
$$\rk (H\cap L) \ge \rk H+\rk L-\rk V \ge m+k-2m \ge 2.
$$
Hence there is an element $\alpha_1\in H\cap L$ such that
$\alpha_1\ne0 $ and $\alpha_1\ne\beta $. If we set
$\alpha_2=\beta+\alpha_1$, then $S_\calu(\alpha_2)$ is the complement
of $S_\calu(\alpha_1)$ relative to $S_\calu(\beta)$. This implies that 
$$\size_\calu(\alpha_1)+\size_\calu(\alpha_2)=\size_\calu(\beta)=k\le2m,$$
so that $\size_\calu(\alpha_j)\le m$ for some $j\in\{1,2\}$. As our
choice of $\alpha_1$ implies that $\alpha_1$ and $\alpha_2$ are both
non-zero, the conclusion of the proposition follows in this case.

{\bf Case II. For every element $\alpha$ of $H$ we have
$\size_\calu(\alpha)\le m+1$.}

If we are in Case II and the conclusion of the proposition does not
hold, then for every $\alpha\in H-\{0\}$ we have $\size_\calu(\alpha)=
m+1$. We shall show this leads to a contradiction.

We consider the collection $\cals=\{S_\calu(\alpha):0\ne \alpha\in H\}$ of
subsets of $\calu$. Each set in $\cals$ has cardinality exactly $m+1$.
If $S$ and $T$ are distinct sets in $\cals$ we may write
$S=S_\calu(\alpha)$ and $T=S_\calu(\beta)$, where $\alpha$ and $\beta$ are distinct
elements of $H-\{0\}$. We then have $\alpha+\beta\in H-\{0\}$, so that
$S_\calu(\alpha+\beta)$ has cardinality $m+1$. But $S_\calu(\alpha+\beta)$ is the
symmetric difference of $S=S_\calu(\alpha)$ and $T=S_\calu(\beta)$, so that
$$\begin{aligned}m+1&=\card{S_\calu(\alpha+\beta)}\\
&=\card{S}+\card{T}-2\card(S\cap T)\\
&=2(m+1)-2\card(S\cap T),
\end{aligned}$$
so that 
\Equation\label{famous scrandy scratch-scratch}\card(S\cap T)=\frac{m+1}2
\EndEquation
for any two distinct sets $S,T\in\cals$.

Since $\rk H=m\ge2$, there exist distinct elements $S$ and $T$ of
$\cals$. It therefore follows from (\ref{famous scrandy
  scratch-scratch}) that $m$ is odd. In particular we have $m\ge3$.

We now apply Fisher's inequality \cite{jukna}*{Theorem 14.6}, which may
be stated as follows. Let $n$ and $k$ be  positive integers, let
$\calu$ be a set of cardinality $n$, and suppose that $\calx$ is a
collection of subsets of $\calu$ such that $\card(S\cap T)=k$ for all
distinct sets $S,T\in\calx$. Then $\card\calx\le n$.

In the present situation, the hypotheses of Fisher's inequality hold
with $n=2m$, $k=(m+1)/2$ and $\calx=\cals$. But if $d\ge m$ is the
dimension of $H$, we have $\card\cals=2^d-1$. Hence Fisher's
inequality gives $2^m-1\le2m$. However, since $m\ge3$ we have
$2m<2^m-1$. This is the required contradiction.
\EndProof

\section{Homology of books of \texorpdfstring{$I$}{I}-bundles}
\label{kalamazoo to timbuctu}

In this section we will use the
  notation introduced in \cite{last}*{Section 2} regarding books of
  $I$-bundles. Recall that if $\calw$ is a book of $I$-bundles then
  $\calb_\calw $ and $\calp_\calw $ denote, respectively the union of all
  bindings of $\calw $ and the union of all its pages; and $|\calw|$
  denotes the manifold $\calb_\calw \cup\calp_\calw $. Each component of
  $\calb_\calw $ is a solid torus.  Each component $P$ of $\calp_\calw $ is
  equipped with the structure of an $I$-bundle over a connected
  $2$-manifold; we denote the associated $\partial I$-bundle by
  $\partial_hP $, and the set
  $\overline{\partial P-\partial_hP}$ by $\partial_vP$ .

\Number\label{torus---not!}
Note that if $F$ is any component of $\partial|\calw|$, then $\chi(F)$ is
the sum of the Euler characteristics of the components of
$\partial_h\calp$ contained in $F$.  Since by \cite{last}*{Definition 2.2} every binding of $\calw$ meets at least one page, $F$
must contain at least one component of $\partial_h\calp$. Hence if
every page of $\calw $ has strictly negative Euler characteristic, then
$\chi(F)<0$ for every component $F$ of $\partial|\calw|$.
\EndNumber

Our next result, Proposition \ref{H_2}, gives a way of computing
$H_2(|\calw|,\calb_0;\Z_2)$, where $\calw$ is a book of $I$-bundles
and $\calb_0$ is a union of certain bindings. For this purpose we need
some notation.

\Notation\label{see one see two} Let $\calw$ be a book of $I$-bundles,
and let $\calb_0$ be a (possibly empty) union of components of
$\calb=\calb_\calw $.  
Let us set $\calb_1=\calb-\calb_0$.  We shall denote by
$C_1(\calw,\calb_0)$ the free $\Z_2$-module generated by the
components of $\calb_1$, and by $C_2(\calw,\calb_0)$ the free
$\Z_2$-module generated by the pages of $\calw$.  

(We think of $C_1(\calw,\calb_0)$ and $C_2(\calw,\calb_0)$ as being
analogous to the groups of $1$-chains and $2$-chains for the chain
complex used to compute the relative homology of a pair of
CW-complexes. Here the bindings and pages of $\calw$ play the roles of
$1$-cells and $2$-cells respectively. From this point of view it is
natural that $C_2(\calw,\calb_0)$ should be independent of
$\calb_0$---as is apparent from the formal definition---since
$|\calb_0|$ contains no pages of $\calw$.)

For each component $B$ of $\calb_1$, let us denote by $d(B)$ the image
in $\Z_2$ of the degree of $B$; and for each component $B$ of
$\calb_1$ and each component $A$ of $\cala_\calw$, let us define
$\delta_{A,B}\in\Z_2$ to be $1$ if $A\subset\partial B$ and $0$
otherwise. We define the boundary homomorphism
$\Delta_{\calw,\calb_0}: C_2(\calw,\calb_0)\to C_1(\calw,\calb_0)$ by
setting $\Delta_{\calw,\calb_0}(P)=\sum_{A,B} \delta_{A,B}d(B) B$ for
each page $P$ of $\calw$, where $A$ ranges over all vertical boundary
annuli of $P$ and $B$ ranges over all components of $\calb_1$.
(Thus the boundary of the $2$-chain $P$ is the formal sum of
  the bindings of odd-valence and odd-degree which are not contained
  in $B_0$ and which meet $P$.)  In the case $\calb_0=\emptyset$ we
shall write $C_1(\calw)$, $C_2(\calw)$ and $\Delta_\calw$ in place of
$C_1(\calw,\emptyset)$, $C_2(\calw,\emptyset)$ and
$\Delta_{\calw,\emptyset}$.
\EndNotation

\Proposition\label{H_2}Suppose that $\calw$ is a book of $I$-bundles,
and that $\calb_0$ is a (possibly empty) union of components of
$\calb=\calb_\calw $. 
Then $H_2(|\calw|,\calb_0;\Z_2)$ is isomorphic to the kernel of
$\Delta_{\calw,\calb_0}: C_2(\calw,\calb_0)\to C_1(\calw,\calb_0)$.
\EndProposition

\Proof
In this proof all homology groups will be understood to have coefficients
in $\Z_2$. We set $C_1=C_1(\calw,\calb_0)$, $C_2=C_2(\calw,\calb_0)$,
and $\Delta=\Delta_{\calw,\calb_0}$. We define $\calb_1$, $d(B)$ and
$\delta_{A,B}$ as in \ref{see one see two}. We set $W=|\calw|$.
Since $H_2(\calb,\calb_0)=0$, we have a natural exact sequence
$$0\longrightarrow H_2(W,\calb_0)\longrightarrow
H_2(W,\calb)\longrightarrow H_1(\calb,\calb_0).$$
Hence $H_2(W,\calb_0)$ is isomorphic to the kernel of the attaching map
$H_2(W,\calb)\to H_1(\calb,\calb_0)$. Setting $\calp=\calp_\calw$ and
$\cala=\cala_\calw$, we have an excision isomorphism
$j:  H_2(\calp,\cala)\to H_2(W,\calb)$, and the domain $H_2(\calp,\cala)$ may be identified with $\bigoplus_P H_2(P,P\cap\cala)$, where $P$ ranges over the pages of $\calw$.
Each summand $H_2(P,P\cap\cala)$ is isomorphic to $H_2(S_P,\partial S_P)$, where $S_P$ denotes the base of the
$I$-bundle $P$, and therefore has rank $1$. If $c_P\in
H_2(W,\calb)$ denotes the image under $j$ of the generator of
$H_2(P,P\cap\cala)$, then the family $(c_P)_P$, indexed by the
pages $P$ of $\calw$, is a basis for $H_2(W,\calb)$. Similarly,
if for each binding $B\subset\calb_1$ we denote by $e_B\in
H_1(\calb,\calb_0)$ the image under inclusion of the generator of the
rank-$1$ vector space $H_1(B)$, then the family $(e_B)_B$, indexed by
the bindings $B$ of $\calb_1$, is a basis for $H_1(\calb,\calb_0)$.
It is straightforward to check that for each
page $P$ of $\calw$, the attaching map $ H_2(W,\calb)\to
H_1(\calb,\calb_0)$ takes $c_P$ to $\sum_{A,B}\delta_{A,B}d(B)e_B$ ,
where $A$ ranges over all vertical boundary annuli of $P$ and $B$
ranges over all components of $\calb_1$. The conclusion of the
Proposition follows.  \EndProof

\Lemma\label{the old concert hall on the bowery}Suppose that $\calw$
is a book of $I$-bundles such that $\chibar(P)=1$ for every page $P$
of $\calw$, and that $\calw_0$ is a connected sub-book of
$\calw $.  Let $p_1$ denote the number of pages of $\calw$ that are
not pages of $\calw_0$. Then $\rk H_1(|\calw|,|\calw_0|;\Z_2)\le2p_1$.
\EndLemma

\Proof
We set $W=|\calw|$ and $W_0=|\calw_0|$.  If we define
$$\chi(W,W_0)=\sum_{i=0}^2(-1)^i\rk H_i(W,W_0;\Z_2)$$
then the exact homology sequence of the pair $(W,W_0)$ implies that
$$\chi(W,W_0)=\chi(W)-\chi(W_0).$$ 
Since each page of $\calw$ has Euler characteristic $-1$, and the
bindings of $\calw$ and their frontiers are of Euler characteristic
$0$, we have $\chi(W)-\chi(W_0)=-p_1$. Hence
$$\rk H_0(W,W_0;\Z_2)-\rk H_1(W,W_0;\Z_2)+\rk H_2(W,W_0;\Z_2)=-p_1.$$

Since $W$ and $W_0$ are connected, we have $H_0(W,W_0;\Z_2)=0$.  To
estimate $\rk H_2(W,W_0;\Z_2)$ we consider the sub-book $\calw_1$ of
$\calw$ consisting of all those pages of $\calw$ that are not
contained in $W_0$, and all bindings of $\calw$ that meet pages not
contained in $W_0$. We set $W_1=|\calw_1|$. We also set
$\calb_0=W_0\cap W_1$, so that $\calb_0$ is in particular the union of a
certain set of bindings of $\calw_1$.  It follows from Proposition
\ref{H_2} that $H_2(W_1,\calb_0;\Z_2)$ is isomorphic to a subspace of
$C_2(\calw,\calb_0)$, the free $\Z_2$-module generated by the pages of
$\calw$.  By definition the dimension of $C_2(\calw,\calb_0)$ is
$p_1$.  Since $H_2(W_1,\calb_0;\Z_2)$ is isomorphic to $H_2(W,W_0;\Z_2)$
by excision, we have $\rk H_2(W,W_0;\Z_2)\le p_1$.
Therefore
$$-p_1=-\rk H_1(W,W_0;\Z_2)+\rk H_2(W,W_0;\Z_2)\le-\rk H_1(W,W_0;\Z_2)+p_1,$$
from which the conclusion of the lemma follows.
\EndProof

\Lemma\label{we'll all go to the seashore}Let $m\ge2$ be an integer.
Suppose that $\calw$ is a book of $I$-bundles such that $\chibar(P)=1$
for every page $P$ of $\calw$, and such that
$\chibar(|\calw|)=2m$. Suppose also that $H_2(|\calw|;\ZZ_2)$ has
dimension at least $m$. Then $\calw$ has a connected sub-book $\calw_0$ such
that $\chibar(|\calw_0|)=m$ and $H_2(|\calw_0|;\Z_2)\ne0$.  \EndLemma

\Proof Since $\chibar(P)=1$
for every page $P$ of $\calw$, the number of pages of every sub-book
$\caly$ of $\calw$ is equal to $\chibar(|\caly|)$. In particular, $\calw$ has
exactly $2m$ pages. In the notation of \ref{see one see two} it
follows that $\rk C_2(\calw)=2m$.

According to Proposition \ref{H_2}, $H_2(|\calw|;\ZZ_2)$ is
isomorphic to the kernel $H$ of $\Delta_{\calw}:
C_2(\calw)\to C_1(\calw)$. The hypothesis of the lemma therefore
implies that $\rk H\ge m$. 

According to the definition of $C_2(\calw)$ (see \ref{see one see
  two}), the set $\calu$ of pages of $\calw$ is canonically identified with a
  basis of $C_2(\calw)$. Since $\rk H\ge m$, 
  Proposition \ref{quality control} gives an
element $\alpha$ of $H$ such that $0<\size_\calu(\alpha)\le m$. 

In the notation of Section \ref{algebra section}, we set
$S=S_\calu( \alpha )$, so that $0<\card{S}\le m$. We define $\calz$
to be the sub-book of $\calw$ whose pages are the elements of $ S $,
and whose bindings are the bindings of $\calw$ that meet pages in the
set $ S $. We set $Z=|\calz|$. We have $Z\ne\emptyset$ since $\card{S}>0$. 

Let $Z_1,\ldots,Z_r$ denote the connected  components  of $Z$,
where $r\ge1$. Then for $i=1,\ldots,r$ we have $Z_i=|\calz_i|$ for
some connected sub-book $\calz_i$ of $\calw$.  We denote by $
S_i\subset S $ the set of all pages of $\calw$ that belong to
$\calz_i$, and we set $ \alpha_i=\sum_{u\in S_i}u $, so that
$\alpha =\alpha _1+\cdots+\alpha _r$.

Let $X$ denote the set of bindings of $\calw$ that are contained in
$\calz$, and for $i=1,\ldots,r$ let $X_i$ denote the set of bindings
of $\calw$ that are contained in $\calz_i$. Let $A$ and $A_i$ denote,
respectively, the subspaces of $C_1(\calw)$ spanned by $X$ and
$X_i$. Then $X$ is the disjoint union of $X_1,\ldots,X_n$ and hence
$A$ is the direct sum of of $A_1,\ldots,A_n$. Since $\alpha \in H$ we
have
$$0=\Delta_{\calw}(\alpha)=\sum_{i=1}^r\Delta_{\calw}(\alpha _i),$$
where $\Delta_{\calw}(\alpha _i)\in A_i$ for $i=1,\ldots,r$. Since the
sum $A_1+\ldots+ A_r$ is direct, it follows that
$\Delta_{\calw}(\alpha _i)=0$ for $i=1,\ldots,r$, so that $\alpha
_1,\ldots,\alpha_r\in H$. The $\alpha_i$ are non-zero since each
component $Z_i$ contains at least one page of $\calw$.

We have $\chibar(Z_1)\le\chibar(Z)=\size_\calu(\alpha)\le m$. 
We set $k=m-\chibar(Z_1)\ge0$, and recursively define sub-books $\caly _j$ of
$\calw$ for $0\le j\le k$, with $\chibar(|\caly _j|)=\chibar(Z_1)+j$, as follows. Set $\caly _0=\calz_1$. If $0\le
j<k$ and
$\caly _j$ has been defined, then $\caly _j$
is a proper sub-book of $\calw$ since
$\chibar(|\caly_j|)=\chibar(Z_1)+j< m<2m=\chibar(W)$. Since $W$ is connected, $|\caly_j|$ must
meet some page
$P$ of $W$ which is not a page of  $\caly_j$. Define $\caly_{j+1}$ to
be the
sub-book of $\calw$ consisting of the pages and bindings of $\caly_j$
together with the page $P$ and the bindings of $\calw$ that meet
$P$.  Then 
$$\chibar(|\caly_{j+1}|)=\chibar(|\caly_j|)+\chibar(P)=\chibar(|\caly_j|)+1=\chibar(Z_1)+j+1,$$
and the recursive definition is complete.

Now set $\calw_0=\caly_{k}$, so that $\chibar(\calw_0|)=m$.
Since the bindings and pages of $\calw_0$ are bindings and pages of
$\calw$, the vector spaces $C_1(\calw_0)$ and $C_2(\calw_0)$ are
naturally identified with subspaces of $C_1(\calw)$ and $C_2(\calw)$.
The boundary homomorphism $\Delta_{\calw_0}:C_2(\calw_0)\to C_1(\calw_0)$ is the
restriction of $\Delta_{\calw}$ to $C_2(\calw_0)$. Hence the kernel of
$\Delta_{\calw_0}$ is $H\cap C_1$. The latter subspace contains the
non-zero element $\alpha_1$, and so $\Delta_{\calw_0}$ has non-trivial
kernel. It now follows from Proposition \ref{H_2} that
$H_2(W_0;\Z_2)\ne0$.
\EndProof

\Lemma\label{all the little children love mary anne}
Suppose
that $\calw$ is a book of $I$-bundles such
that $\chibar(P)>0$ for every page of $\calw$.
Then there  is a book of $I$-bundles $\calw'$ such
that 
that $|\calw'|=|\calw|$, and such that $\chibar(P)=1$ for every page
$P$ of $\calw$.
\EndLemma

\Proof Set $W=|\calw|$ and $\calp=\calp_\calw$. Let $S$ denote the
base of the $I$-bundle $\calp$, and let $q:\calp\to S$ denote the
bundle map. Since every component of $S$ has negative Euler
characteristic, there is a closed $1$-manifold $\calc\subset S$ such
that every component of $S-\calc$ has Euler characteristic $-1$. Let
$\caln$ be a regular neighborhood of $\calc$ in $S$. Set
$\calb'=q^{-1}(\caln)$ and $\calp'=q^{-1}(\overline{S-\caln})$. Then
$\calp'$ inherits an $I$-bundle structure from $\calp$, and we need
only set $\calw'=(W,\calb',\calp')$.  \EndProof

\Lemma\label{running in and out the hut}
Let $m\ge1$ be an integer. Suppose that $M$ is a closed, orientable,
irreducible $3$-manifold which is $(m+1)$-small (\ref{who the hell
  said anything about christmas}).  Suppose that $\calw$ is a book of
$I$-bundles with $W\doteq|\calw|\subset M$, that $\chibar(P)=1$ for
every page $P$ of $\calw$, and that $\chibar(|\calw|)=2m$. Suppose
also that $H_2(|\calw|;\ZZ_2)$ has dimension at least $m$. Then
$\calw$ has a sub-book $\calw_0$ such that \Conclusions
\item $\chibar(|\calw_0|)=m$, and
\item the inclusion homomorphism $H_1(|\calw_0|;\Z_2)\to H_1(M;\Z_2)$
  is {\it either} surjective {\it or}
 has image of rank at most $\max(m,2)$.
\EndConclusions
\EndLemma

\Proof We first consider the case $m\ge2$. In this case, according to
Lemma \ref{we'll all go to the seashore}, $\calw$ has a connected
sub-book $\calw_0$ such that $W_0=|\calw_0|$ satisfies $\chibar(W_0)=m$ and
$H_2(W_0;\Z_2)\ne0$. 

If it happens that the inclusion homomorphism $\pi_1(W_0)\to \pi_1(M)$
is surjective, then in particular the inclusion homomorphism
$H_1(W_0;\Z_2)\to H_1(M;\Z_2)$ is surjective, so that the conclusion
of the lemma holds.

Now suppose the inclusion homomorphism $\pi_1(W_0)\to \pi_1(M)$
  is not surjective. According to 
\ref{torus---not!}, no component of $\partial W_0$ is a sphere. Hence,
in the notation of \ref{all chalk, no traction} we have $W_0\in\calx_M$.

The manifold $W_0$ is not a handlebody, since $H_2(W_0;\Z_2)\ne0$.
According to \ref{torus---not!}, no component of $\partial W_0$
is a torus. The hypotheses of Proposition \ref{newest old prop 2} are
now seen to hold with $Y=W_0$ and $c=m$. (The condition that $M$ is $
(m+1)$-small is a hypothesis of the present lemma.)  It therefore
follows from Proposition \ref{newest old prop 2} that $i(W_0)\ge0$.
According to the definition of $i(X)$ given in \ref{all chalk, no
  traction}, this means that $r(W_0)\le\chibar(W_0)$, where $r(W_0)$
is the rank of the inclusion homomorphism $\pi_1(W_0)\to \pi_1(M)$. In
particular, the inclusion homomorphism $H_1(W_0;\Z_2)\to H_1(M;\Z_2)$
has rank at most $\chibar(W_0)=m$. This completes the proof of the
lemma in the case $m\ge2$.

We now consider the case $m=1$. In this case we select a page $P_0$ of
$\calw$ and define the sub-book $\calw_0$ to consist of $P_0$ and the
bindings of $\calw$ that meet $P_0$. Then $W_0\doteq|\calw_0|$ is
connected and $\chibar(W_0)=1$.

If it happens that the inclusion homomorphism $\pi_1(W_0)\to \pi_1(M)$
is surjective, then in particular the inclusion homomorphism
$H_1(W_0;\Z_2)\to H_1(M;\Z_2)$ is surjective, so that the conclusion
of the lemma holds.

Now suppose the inclusion homomorphism $\pi_1(W_0)\to \pi_1(M)$ is not
surjective. According to \ref{torus---not!}, no component of
$\partial W_0$ is a sphere. Hence, we have $W_0\in\calx_M$.  The
hypotheses of Proposition \ref{new but not newest old prop 2} are now
seen to hold with $Y=W_0$ and $c=1$.  It therefore follows from
Proposition \ref{new but not newest old prop 2} that $i(W_0)\ge-1$,
i.e. that $r(W_0)\le\chibar(W_0)+1=2$. In particular, the inclusion
homomorphism $H_1(W_0;\Z_2)\to H_1(M;\Z_2)$ has rank at most $2$. Thus
the lemma is proved in all cases.
\EndProof

\Proposition\label{sifting sand}
Let $m\ge1$ be an integer. Suppose that $M$ is a closed, orientable,
irreducible $3$-manifold which is $(m+1)$-small.  Suppose that
$\calw$ is a book of $I$-bundles with $W\doteq|\calw|\subset M$, that
$\chibar(P)>0$ for every page $P$ of $\calw$, and that
$\chibar(|\calw|)=2m$. Then the image of the inclusion homomorphism
$H_1(W;\Z_2)\to H_1(M;\Z_2)$ has dimension at most $\max(3m,4)$.
\EndProposition

\Proof
According to Lemma \ref{all the little children love mary anne} we may
assume without loss of generality that $\chibar(P)=1$ for every page
$P$ of $\calw$.

We shall let $T$ denote the image of the inclusion homomorphism
$j:H_1(W;\Z_2)\to H_1(M;\Z_2)$.

We consider first the case in which $H_2(W ;\ZZ_2)$ has dimension
at most $m-1$. In this case we note that 
$$2m=\chibar(W)
=-\rk H_0(W;\Z_2)+\rk H_1(W;\Z_2)-\rk H_2(W;\Z_2)
\ge -m + \rk W,$$ so
that $\rk W\le3m$. It follows immediately that $\rk T\le3m$ in this
case.

There remains the case in which $H_2( W;\ZZ_2)$ has dimension at
least $m$. In this case, according to Lemma \ref{running in and out
the hut}, there is a sub-book $\calw_0$ of $\calw$ such that
$\chibar(|\calw_0|)=m$, and such that the inclusion homomorphism $j_0:H_1(|\calw_0|;\Z_2)\to
H_1(M;\Z_2)$ either 
is surjective or has image of rank at most $\max(m,2)$.

Set $W_0 = |\calw_0|$.  By \cite{last}*{Lemma 2.11}, we have
$$\rk( W_0 )\le2\barchi( W)+1=2m+1.$$
Hence in the subcase where $j_0$ is surjective, we have $\rk
H_1(M;\Z_2)\le \rk H_1(W_0;\Z_2)\le2m+1\le 3m$; since $T$ is a
subspace of $H_1(M;\Z_2)$, we in particular have $\rk T\le 3m$ in this
subcase.

Finally we consider the subcase in which $T_0\doteq j_0(H_1( W_0
;\Z_2))$ has dimension at most $\max(m,2)$. Since $\chibar(W)=2m$ and
$\chibar(W_0)=m$, and since $\chibar(P)=1$ for each page of $W$, the
number of pages of $\calw$ that are not pages of $\calw_0$ is equal to
$m$. Hence by Lemma \ref{the old concert hall on the bowery}, we have
$\rk H_1(W,W_0;\Z_2)\le 2m$.

Let $L$ denote the cokernel of the inclusion homomorphism
$H_1(W_0;\Z_2)\to H_1(W;\Z_2)$.  The natural surjection from
$H_1(W;\Z_2)$ to $T$ induces a surjection from $L$ to $T/T_0$. Hence
$$\rk T-\rk T_0=\rk(T/T_0)\le\rk L\le\rk
H_1(W,W_0;\Z_2)\le2m.$$
Since $\rk T_0\le\max(m,2)$, it follows that
$$\rk T\le2m+\max(m,2)=\max(3m,4),$$
as required.
\EndProof

\section{De-singularizing surfaces}\label{illegitimis non carborundum}

This section is devoted to the proof of Theorem \ref{top 6}, which was
stated in the introduction.

\Proof[Proof of Theorem \ref{top 6}]
We use the terminology of \cite{last}.  Applying
\cite{last}*{Proposition \towerproposition}, we find a good tower
$${\mathcal T}=(M_0,N_0,p_1,M_1,N_1,p_2,\ldots,p_n,M_n,N_n),$$ with
base $M_0$ homeomorphic to $M$ and with some height $n\ge0$, such that
{$N_n$} contains a connected incompressible closed surface $F$ of genus $\le
g$.  According to the definition of a good tower, $\partial N_{n}$ is
incompressible (and, {\it a priori}, possibly empty) in $M_{n}$. Hence
$N_{n}$ is $\pi_1$-injective in $ M_{n}$. Since $F$ is incompressible
in $N_{n}$, it follows that it is also incompressible in $M_{n}$.

Since $M$ is simple it follows from \cite{last}*{Lemma \simpletower}
that all the $M_j$ and $N_j$ are simple.

Let $k$ denote the least integer in $\{0,\ldots,n\}$ for which $M_k$
contains a closed incompressible surface $S_k$  of genus at most
$g$.  To prove the theorem it suffices to show that $k=0$.
Let $h$ denote the genus of $S_k$. Since $M_k$ is simple we have $h\ge2$.

Suppose that $k\ge1$. The minimality of $k$ implies that $M_{k-1}$
contains no closed incompressible surface of genus at most $g$. In
particular:

\Claim\label{will it core a apple} $M_{k-1}$ contains no
closed incompressible surface of genus at most $h$.
\EndClaim

From \ref{will it core a apple} it follows that, in particular,

\Claim\label{tall skinny cappucino}$M_{k-1}$
is $h$-small. 
\EndClaim

We now evoke \cite{last}*{Proposition \newimprovedoldpropthree}, which
states that if $\tN$ is a $2$-sheeted covering of a simple, compact,
orientable $3$-manifold $N$, and if $ \tN $ contains a closed,
incompressible surface of a given genus $h\ge2$, then either (1) $N$
contains a closed, connected, incompressible surface of genus at most
$h$, or (2) $N$ is closed and there is a connected book of
$I$-bundles $\calw$ with $W=|\calw|\subset N $ such that
$\chibar(W)=2h-2$, every page of $\calw$ has strictly negative Euler
characteristic, and every component of $\overline{N-W}$ is a
handlebody.   Observe that the hypotheses of
\cite{last}*{Proposition \newimprovedoldpropthree} hold in the present
situation if we set $N=N_{k-1}$ and $\tN=M_k$.

Suppose that alternative (1) of the conclusion of
\cite{last}*{Proposition \newimprovedoldpropthree} holds in the
present situation, i.e.  that $ N_{k-1}$ contains an incompressible
closed surface $S_{k-1}$ with $\genus(S_{k-1})\le h\le g$.  According
to the definition of a good tower, $\partial N_{k-1}$ is an
incompressible (and possibly empty) surface in $M_{k-1}$. Hence
$N_{k-1}$ is $\pi_1$-injective in $ M_{k-1}$. Since $S_{k-1}$ is
incompressible in $N_{k-1}$, it follows that it is also incompressible
in $M_{k-1}$. This contradicts \ref{will it core a apple}.

Now suppose that alternative (2) of the conclusion
  of \cite{last}*{Proposition
  \newimprovedoldpropthree} holds in the present situation, i.e.:
\Claim\label{yes it will core a apple} $N_{k-1}$ is closed and there
is a connected book of $I$-bundles $\calw$ with $W=|\calw|\subset
N_{k-1} $ such that $\chibar(W)=2(h-1)$, every page of $\calw$ has
strictly negative Euler characteristic, and every component of
$\overline{N_{k-1}-W}$ is a handlebody.
\EndClaim

Since $N_{k-1}$ is closed we have $N_{k-1}=M_{k-1}$.

It now follows from \ref{tall skinny cappucino} and \ref{yes it will
  core a apple} that the hypotheses of Proposition \ref{sifting sand}
hold with $m=h-1$, and with $M_{k-1}$ in place of $M$. Hence by
Proposition \ref{sifting sand}, the image of the inclusion
homomorphism $H_1(W;\Z_2)\to H_1(M_{k-1};\Z_2)$ has dimension at most
$\max(3h-3,4)$.  On the other hand, since by \ref{yes it will core a
  apple} every component of $\overline{N_{k-1}-W}$ is a handlebody,
the inclusion homomorphism $H_1(W;\Z_2)\to H_1(M_{k-1};\Z_2)$ is
surjective. Hence
 $$\rk M_{k-1} \le\max(3h-3,4) \le \max(3g-3,4).$$
On the other hand, since by hypothesis we have $\rk
M_0\ge\max(3g-1,6)$, it follows from \cite{last}*{Lemma \tworminusfour}
that for any index $j$ such that $0\le j\le n$ and such that $M_j$ is
closed, we have $\rk M_j\ge\max(3g-2,5)$. This is a contradiction, and
the proof is complete.
\EndProof

\section{An example}\label{example section}

In this section we investigate the extent to which Theorem \ref{top 6}
is sharp. Our discussion focuses on the case $g=2$ of Theorem \ref{top
  6}, although the methods can be applied more generally.  To show
that the theorem is sharp for $g=2$ one would need an example of a
closed simple $3$-manifold $M$ with $\rk M = \max(3g-2, 5) = 5$, such
that $\pi_1(M)$ contains a genus-$2$ surface group but $M$ contains no
closed, incompressible surface of genus $2$. Proposition \ref{example}
below asserts the existence (and the proof gives an explicit example)
of a closed simple $3$-manifold $M$ with $\rk M =4$, such that
$\pi_1(M)$ contains a genus-$2$ surface group but $M$ contains no
closed, incompressible surface whatever. We will also show why our
construction cannot give a similar example in which $\rk M$ is $5$
rather than $4$; however, we have no reason to think that such an
example does not exist.

Our example is based on a Dehn surgery construction, and we shall use
notation that is standard in the study of Dehn surgery. If $Q$ is a
compact, orientable $3$-manifold whose boundary is a torus, we define
a {\it slope} for $Q$ to be an isotopy class of unoriented simple
closed curves in $Q$. If $\alpha$ and $\beta$ are slopes, we denote
their geometric intersection number by $\Delta(\alpha,\beta)$. We
define an {\it essential surface} in $Q$ to be a $\pi_1$-injective,
properly embedded, orientable surface which is not
boundary-parallel. If $S$ is an essential surface, all its boundary
components represent the same slope, called the {\it boundary slope}
of $S$.

The following result is essentially due to Cooper-Long and Li.

\Theorem\label{li} Let $Q$ be a simple $3$-manifold whose boundary is
a single torus. Let $S\subset Q$ be an essential
surface with two boundary components. Suppose that $S$ is
not a fiber or semifiber. Let $s$ denote its boundary
slope. Then there is an integer $\Gamma$ such that for every slope $r$
with $\Delta(r,s)\ge \Gamma$, the fundamental group of the Dehn-filled
manifold $Q(r)$ contains an isomorphic copy of $\pi_1(T)$, where $T$
is a closed orientable surface with $\chi(T)=\chi(S)$.
\EndTheorem

\Proof This follows from the proof of \cite{li}*{Theorem 1.2}. (See
also \cite{cooperlong1} and \cite{cooperlong2}.) The statement of
\cite{li}*{Theorem 1.2} does not contain the information that
$\chi(T)=\chi(S)$, but it follows from the proof because $T$ is
constructed from $S$, as in \cite{freedmanfreed}, by adding a singular
annulus joining the two boundary components of $S$.
\EndProof

Theorem \ref{li} will be applied via the following result:

\Proposition\label{li-figue, li-raisin} Let $Q$ be a simple
$3$-manifold whose boundary is a single torus. Suppose that $Q$
contains no closed incompressible surface of genus $>1$. Let $S\subset
Q$ be a separating essential surface with
$\chi(S)=-2$. Suppose that $S$ is not a semifiber. Then
\begin{enumerate} 
\item\label{no more nor four} $Q$ has Heegaard genus $\le 4$; and
\item\label{far more} there are infinitely many
slopes $r$ such that $M:=Q(r)$ has the following properties:
\begin{itemize}
\item $\pi_1(M)$ contains a genus-$2$ surface group;
\item $M$ contains no closed incompressible surface; and
\item $\rk M=\rk Q$.
\end{itemize}
\end{enumerate}
\end{proposition}

\Proof To prove that Conclusion (\ref{no more nor four}) holds we will
construct a Heegaard splitting of the form $Q = V\cup W$ where $V$ is
a compression body and $W$ is a handlebody of genus $4$.

Let $A$ denote the union of three disjoint properly embedded
arcs in $S$ such that $S-A$ is simply-connected.  Let $V$ be a regular
neighborhood of $\partial M \cup A$.  Then $V$ is a compression
body such that $\partial_-V = \partial M$ and $\partial_+V$ has genus
$4$. By adding a $2$-handle to $V$ one obtains a regular neighborhood
$N$ of $\partial M \cup S$.  The frontier of $N$ consists of two
surfaces $F_1$ and $F_2$ of genus $2$.  Since $Q$ is simple and contains no
incompressible surface of genus $>1$, $Q-N$ is a union of two disjoint
handlebodies of genus $2$.  Thus $W = Q-V$ consists of two
handlebodies joined by a $1$-handle, and hence is a handlebody of
genus $4$.

To prove (\ref{far more}), we let $s$ denote the boundary slope of
$S$.  Let $\alpha$ be an indivisible element of $H_1(\partial
M;\ZZ_2)$ which belongs to the kernel of the inclusion homomorphism
$H_1(\partial M;\ZZ_2)\to H_1(M;\ZZ_2)$. Let us extend
$\alpha$ to a basis $(\alpha,\beta)$ for $H_1(\partial M;\ZZ_2)$. As
there are infinitely many choices for $\beta$ we may take $\beta\ne
s$. For each positive integer $n$, the primitive homology class
$\alpha+2n\beta$ determines a slope $r_n$. Since $\alpha+2n\beta$
lies in the kernel of the inclusion homomorphism $H_1(\partial
M;\ZZ_2)\to H_1(M;\ZZ_2)$, we have $\rk Q(r_n)=\rk Q$ for
each $n$. On the other hand, we have
$$\Delta(r_n,s)\ge n\Delta(\beta,s)-\Delta(\alpha,s)$$
for each $n$. Here $\Delta(\beta,s)\ne0$ since $\beta\ne s$, and so
$\Delta(r_n,s)\to\infty$ as $n\to\infty$. Hence Theorem \ref{li}
guarantees that for any sufficiently large $n$ the group
$\pi_1(Q(r_n))$ contains an isomorphic copy of $\pi_1(T)$, where $T$
is a closed orientable surface with $\chi(T)=\chi(S)=-2$; that is,
$\pi_1(Q(r_n))$ contains a genus-$2$ surface group for all
sufficiently large $n$.  

On the other hand, by a theorem of Hatcher \cite{hatcher}, there are
only finitely many boundary slopes for $M$. Since $Q$ is simple and
contains no closed incompressible surface of genus $>1$, the manifold
$Q(r)$ cannot contain a closed incompressible surface unless $r$ is a
boundary slope. Hence for sufficiently large $n$ the manifold
$Q(r_n)$ contains no closed incompressible surface.
\EndProof

The next result produces our example.

\Proposition\label{example} There exists a simple, closed, orientable
$3$-manifold $M$ with $\rk M = 4$, such that $\pi_1(M)$ contains a
genus-$2$ surface group but $M$ contains no incompressible surface.
\EndProposition

\Proof
The Hodgson-Weeks census of cusped hyperbolic $3$-manifolds has been
extended by Thistlethwaite \cite{morwen} to include manifolds which have ideal
triangulations with eight tetrahedra.  We let $\Theta$ denote
the ideal-triangulated manifold {\tt t12045} in the Thistlethwaite
census. 

The program SnapPy
\cite{snappy} reports that $\Theta$ is hyperbolic with
finite volume and one cusp, and that $H_1(Q;\ZZ)$ is isomorphic to
$\ZZ_2\oplus\ZZ_2\oplus\ZZ_2\oplus\ZZ$. Using the program t3m \cite{t3m} to enumerate spunnormal surfaces, in the sense
of \cite{walsh}, with respect to $\calt$, one finds
a surface $\Sigma_0$ with Euler characteristic $-1$ and one end. 

The compact core $Q$ of $\Theta$ is a simple manifold with one
boundary torus.  Truncating $\Sigma_0$ gives a properly embedded
surface $S_0\subset Q$ having Euler characteristic $-1$ and one
boundary component.  Dehn filling on the boundary slope of $S_0$
produces a manifold with first homology
$\ZZ_2\oplus\ZZ_2\oplus\ZZ_2\oplus\ZZ_2$.  In particular, the boundary
curve of $S_0$ is non-trivial in $H_1(Q;\ZZ)$, so $S_0$ is a Klein
bottle with one disk removed.  We let $S$ denote the frontier of a
regular neighborhood $V$ of $S_0$, so that $S$ is an orientable
surface with two boundary components and genus $1$.

The t3m program reports that Thistlethwaite's triangulation $\calt$ of
$\Theta$ admits a taut structure, in the sense of \cite{lackenby}.
The definition of a taut structure involves an assignment of a
transverse orientation to every $2$-simplex of $Q$. One of the
conditions that these transverse orientations are stipulated to
satisfy is that every $3$-simplex has two faces for which the
transverse orientation is inward and two for which it is outward. In
particular each $3$-simplex has a distinguished pair of opposite
edges, namely the common edge of the two outward faces and the common
edge of the two inward faces. Thus there is a distinguished normal
quadrilateral type in each $3$-simplex, namely the one which is
disjoint from the distinguished edges.

The t3m program verifies that for a suitable taut structure on
$\calt$, the spunnormal surface $\Sigma_0$ has the property that all
of its quadrilaterals are of distinguished type. It is clear that $S$
may be obtained by truncating a spunnormal surface $\Sigma$ which has
the same quadrilateral types as $\Sigma_0$. In particular all the
quadrilaterals of $\Sigma$ are of distinguished type.

An unpublished theorem of Dunfield's \cite{dunfield} implies that if
an orientable spunnormal surface in a taut ideal triangulation has the
property that all its quadrilaterals are of distinguished type, then
the properly embedded surface obtained from it by truncation is
essential.  Thus we see that $S$ is essential.

The surface $S$ separates $Q$ since it is the frontier of $V$. If $S$
is a semifiber then $W:=\overline{Q-V}$ is a twisted $I$-bundle over a
surface with associated $\partial I$-bundle $S$, and hence
$H_1(W,S;\ZZ_2)\cong\ZZ_2$. By excision it follows that
$H_1(Q,V;\ZZ_2)\cong\ZZ_2$. Since $Q$ and $V$ are connected, it
follows from the long exact homology sequence of the pair $(Q,V)$ that
$\rk Q\le1+\rk V=1+\rk S_0=3$. This is a contradiction since we have
seen that $\rk Q=4$. Thus we have shown that $S$ is not a semifiber.

The t3m program also verifies that all closed spunnormal surfaces with
respect to the ideal triangulation $\mathcal T$ bound handlebodies,
and hence are compressible.  Hence $Q$ has no closed incompressible
surfaces.

It now follows from Proposition \ref{li-figue, li-raisin} that there
are infinitely many distinct Dehn surgeries on $Q$ which produce
manifolds $M$ with $\rk M=4$, such that $\pi_1(M)$ contains a
genus-$2$ surface group but $M$ contains no closed incompressible
surface.
\EndProof

\Remark The proof of Proposition \ref{example} that we have given
requires constructing the manifold $M$ by a Dehn filling from a
manifold $Q$ satisfying the hypotheses of Proposition \ref{li-figue,
  li-raisin}. Conclusion (\ref{no more nor four}) of Proposition
\ref{li-figue, li-raisin} asserts that any such manifold $Q$ must have
a Heegaard splitting of genus at most $4$.  Since $Q$ has connected
boundary, one of the two compression bodies in this splitting will be
a handlebody.  Thus $Q$ is obtained from a handlebody of genus at most
$4$ by adding $2$-handles. This implies that $\rk Q\le4$, and hence
that $\rk M\le4$ for any manifold $M$ obtained from $Q$ by a Dehn
filling. This is why our method cannot furnish an example with $\rk
M=5$, as it would have to do in order to show that Theorem \ref{top 6}
is sharp when $g=2$.
\EndRemark

\Remark Thurston's Dehn filling theorem implies that the proof of
Proposition \ref{example} gives infinitely many non-homeomorphic
manifolds with the stated properties.
\EndRemark

\section{Non-fibroid surfaces}\label{non-fibroid section}

In this section we will establish a slightly stronger version of
Theorem \ref{top 6}, Proposition \ref{stealing that extra bow}, which
will be useful for volume estimates.

\Definition
Following the terminology that we introduced in \cite{nonsep}, we
define a {\it fibroid} in a closed, orientable topological 3-manifold
$M$ to be a closed incompressible surface $S\subset M$ such that each
component of the manifold-with-boundary obtained by splitting $M$
along $S$ has the form $|\calw|$ for some book of $I$-bundles $\calw$
whose pages are all of negative Euler characteristic.
\EndDefinition

\Proposition\label{stealing that extra bow}
Let $g$ be an integer $\ge2$.  Let $M$ be a closed simple $3$-manifold
such that $\rk M \ge \max(3g-1,6)$ and $\pi_1(M)$ has a subgroup
isomorphic to a genus-$g$ surface group. Then $M$ contains a closed,
incompressible surface which has genus at most $g$ and is not a
fibroid.
\EndProposition

\Proof
According to Theorem \ref{top 6}, $M$ contains a closed,
incompressible surface of some genus $h\le g$. We may take $h$ to be
minimal in the sense that $M$ contains no closed, incompressible
surface of genus $<h$. Since $M$ is simple we have $h\ge2$. We
distinguish two cases.

{\bf Case I. There is a separating closed incompressible surface
  $S\subset M$ with genus $h$.}  We shall show that $S$ is not a
fibroid.  Let $W$ and $W'$ denote the closures of the components of
$M-S$.  We have $\chibar(W)=\chibar(W')=h-1$. Suppose that $F$ is a
fibroid, so that there are books of $I$-bundles $\calw$ and $\calw'$
whose pages are all of negative Euler characteristic, such that
$|\calw|=W$ and $|\calw'|=W'$. It then follows from
\cite{last}*{\moosday} that $\rk W\le 2\chi(W)+1=2h-1$
and $\rk W'\le 2\chi(W')+1=2h-1$.

Consider the Mayer-Vietoris exact sequence
$$
\xymatrix@-4pt{
H_1(F)
\ar[rr]^-{\iota_*+\iota'_*} &&
H_1(W)\oplus H_1(W')
\ar[r]^-{\alpha} &
H_1(M)
\ar[r]^\beta &
H_0(F)
\ar[rr]^-{\iota_*+\iota'_*} &&
H_0(W)\oplus H_0(W')
}
$$ 
where coefficients are taken in $\Z_2$, and where $\iota$ and $\iota'$
are the inclusions of $F$ into $W$ and $W'$.  Since $F$ and $W$ are
path-connected, $\iota_*:H_0(F)\to H_0(W)$ is an isomorphism; hence
$\beta=0$, and $\alpha$ is surjective. It is a standard consequence of
Poincar\'e-Lefschetz duality that the dimension of
$\iota_*(H_1(F))\subset H_1(W)$ is equal to the genus $h$ of
$F=\partial W$.  Hence $(\iota_*+\iota'_*)(H_1(F))$
is a subspace of dimension at least $h$ in
$H_1(W)\oplus H_1(W')$. It follows that
$$\begin{aligned}
\rk M
&\le\rk(H_1(W;\Z_2)\oplus H_1(W';\Z_2)) - h \cr
&\le(2h-1)+(2h-1)-h\cr
&=3h-2\cr
&\le3g-2,
\end{aligned}$$
which contradicts the hypothesis.

{\bf Case II. There is no separating closed incompressible surface of
genus $h$ in $M$. } By our choice of $h$, there is also no 
closed incompressible surface of genus $<h$ in $M$. By definition this
means that $M$ is $h$-small. 

Our choice of $h$ also guarantees that there is a
closed incompressible surface $S$ of genus $h$ in $M$. Since we are in
Case II, the surface $S$ is non-separating. We shall show that $S$ is
not a fibroid.

Fix a regular neighborhood $N$ of $S$ in $M$, and set
  $W=\overline{M-N}$. Since $S$ is non-separating,  $W$ is connected. We
  have $\chibar(W)=2h-2$. Suppose that $S$ is a fibroid,
  so that there is a book of $I$-bundles
$\calw$ whose pages are all of negative Euler
characteristic, such that
 $|\calw|=W$. Since $M$ is $h$-small, the hypotheses of Proposition
  \ref{sifting sand} are now seen to hold with $m=h-1$. Hence  if $T$
  denotes the
  image of the inclusion homomorphism
$H_1(W;\Z_2)\to H_1(M;\Z_2)$, it follows from Proposition
  \ref{sifting sand} that $T$ has dimension at most $\max(3h-3,4)$.

If $c$ is the class in $H_1(M;\Z_2)$ defined by a simple closed curve
that crosses $S$ in one point, then $H_1(M;\Z_2)$ is spanned by $c$
and $T$. It follows that $H_1(M;\Z_2)$ has dimension at most
$\max(3h-2,5)$. This contradicts the hypothesis.
\EndProof

\section{Volumes}\label{volume section}

In this section we will establish Theorem \ref{geom 6} which was
stated in the introduction. One of the ingredients is a result due to
Agol, Storm, and Thurston from \cite{ast}. The information from
\cite{ast} that we need is summarized in Theorem \firstAST\ of
\cite{last}, which can be paraphrased as saying that if $M$ is a
closed orientable hyperbolic $3$-manifold containing a connected
incompressible closed surface which is not a fibroid, then
$\vol(M)>3.66$.

We also recall that a group $\Gamma$ is said to be $k$-free, where $k$
is a positive integer, if every subgroup of $\Gamma$ having rank at
most $k$ is a free group. The following result provides the
transition between the earlier sections of this paper and the
applications to volumes, which include the proofs of Theorem \ref{geom
  6} and of the corresponding result in \cite{4-free}.

\Proposition\label{if not why not}
Let $k\ge3$ be an integer, and let $M$ be a closed orientable simple
$3$-manifold such that $\rk M\ge\max(3k-4,6)$. Then either $\pi_1(M)$
  is $k$-free, or $M$ contains a closed incompressible surface of
  genus at most $k-1$ which is not a fibroid.
\EndProposition

\Proof First consider the case in which $\pi_1(M)$ has a subgroup
isomorphic to a genus-$g$ surface group for some $g$ with $1<g\le k-1$. The
hypothesis then implies that $\rk M\ge\max(3g-1,6)$, and it follows
from Proposition \ref{stealing that extra bow} that $M$ contains a
closed, incompressible surface which is not a fibroid and has genus at
most $g\le k-1$.

Now consider the case in which $\pi_1(M)$ has no subgroup isomorphic
to a genus-$g$ surface group for any $g$ with $1<g\le k-1$.  In this
case, since $\rk M \ge k+2$, it follows from \cite{accs}*{Proposition
  7.4 and Remark 7.5} that $\pi_1(M)$ is $k$-free.
\EndProof

We now turn to the proof of Theorem \ref{geom 6}.

\Proof[Proof of Theorem \ref{geom 6}]
 Assume that $\rk M\ge6$. Then according to Proposition \ref{if
not why not}, either $\pi_1(M)$ is $3$-free, or $M$ contains a closed
incompressible surface of genus at most $2$ which is not a fibroid. If
$\pi_1(M)$ is $3$-free, it follows from \cite{last}*{Corollary
\threefreevolume} that $\vol(M)>3.08$. If $M$ contains a closed
incompressible surface which is not a fibroid, it
follows from \cite{last}*{Theorem \firstAST} that $\vol(M)>3.66$. In
either case the hypothesis is contradicted.
\EndProof

\end{document}